# RESEARCH ANNOUNCEMENT

# SECONDARY INVARIANTS AND THE SINGULARITY OF THE RUELLE ZETA-FUNCTION IN THE CENTRAL CRITICAL POINT

ANDREAS JUHL

ABSTRACT. The Ruelle zeta-function of the geodesic flow on the sphere bundle $S(X)$ of an even-dimensional compact locally symmetric space $X$ of rank 1 is a meromorphic function in the complex plane that satisfies a functional equation relating its values in $s$ and $-s$. The multiplicity of its singularity in the central critical point $s = 0$ only depends on the hyperbolic structure of the flow and can be calculated by integrating a secondary characteristic class canonically associated to the flow-invariant foliations of $S(X)$ for which a representing differential form is given.

Let $Y = G/K$ be a rank one symmetric space of the non-compact type, i.e., $Y$ is a real, complex, or quaternionic hyperbolic space or the (16-dimensional) hyperbolic Cayley-plane. Let $\Gamma$ be a uniform lattice in the (connected simple) isometry group $G$ of $Y$ without torsion. $\Gamma$ acts properly discontinuous on $Y = G/K$ ($K$ a maximal compact subgroup of $G$), and $X = \Gamma\backslash G/K$ is a compact locally symmetric space. We consider $X$ as a Riemannian manifold with respect to an arbitrary (constant) multiple $g$ of the metric $g_0$ induced by the Killing form on the Lie algebra of $G$. Then the Riemannian manifold $(X, g)$ is a space of negative curvature.

The negativity of the curvature of the metric $g$ on $X$ implies the existence of an infinite countable set of prime closed geodesics in $X$ with a discrete set of prime periods accumulating at infinity.

Let $\Phi_t$ be the geodesic flow on the unit sphere bundle $S(X)$ of the space $(X, g)$.

The prime period of a periodic orbit of $\Phi_t$ on $S(X)$ coincides with the length of the closed geodesic in $X$ obtained by projecting the periodic orbit into $X$. Now we use these periods to define the zeta-function

$$Z_R(s) = \prod_c (1 - \exp(-sl_c))^{-1} \tag{1}$$

for $s \in \mathbb{C}$ such that $\mathrm{Re}(s) > h$ ($h$ being the topological entropy of the geodesic flow $\Phi_t$ on $S(X)$). The product in (1) runs over all closed *oriented* geodesics $c$ in $X$, and $l_c$ denotes the length of $c$ as a curve in $X$. Note that for each (unoriented) closed geodesic $c$ in $X$ there are two lifts of $c$ as periodic orbits of $\Phi_t$ in correspondence with the two possibilities to orient $c$.

The function $Z_R$ is well known as the Ruelle zeta-function of the geodesic flow [5]. The Euler product (1) defines a holomorphic function in the half-plane $\mathrm{Re}(s) > h$.

Now by symbolic dynamics the zeta-function $Z_R$ can be written for large $\mathrm{Re}(s)$ as an alternating product of zeta-functions associated to suspensions of subshifts of finite type. The latter zeta-functions coincide with Fredholm-determinants $\det(1 - L_*(s))$ (in the sense of Grothendieck) of holomorphic families $L_*(s)$ of nuclear transfer operators on certain spaces of differential forms. This representation implies that $Z_R$ has a meromorphic continuation to the complex plane (see [5], [14]). While these arguments establish the existence of a meromorphic continuation, it seems to be rather difficult to prove results on the positions and the multiplicities of the singularities of $Z_R$ by the same method.

On the other hand, the zeta-function $Z_R$ can be written as a product of generalized Selberg zeta-functions. Generalized Selberg zeta-functions are also defined by Euler products similar to (1) but with more complicated local Euler factors containing monodromy contributions of the loops $c$ in certain vector bundles on $S(X)$ (see [3], [5], [7], [12], [17]). The generalized Selberg zeta-functions in turn can be investigated by using trace formula techniques which are at present the only known methods to uncover the deeper relations between periodic orbits and the geometry and topology of the underlying space. However, in contrast to the dynamical point of view, along the usual trace formula arguments the relation between the hyperbolic structure and the analytical properties of the zeta-function remains mysterious.

The main result of the present note describes a direct relation between the multiplicity of the singularity of $Z_R$ in the central critical point $s = 0$ and the hyperbolic structure of the flow.

Let us first consider the well-known special case of the Ruelle zeta-function of the geodesic flow of a compact Riemannian surface $X = \Gamma \backslash H^2$ of negative Euler characteristic $\chi(X)$. Consider $X$ as a Riemannian manifold with the metric inherited from the hyperbolic metric

$$y^{-2}(dx^2 + dy^2) \tag{2}$$

of constant curvature $-1$ on the upper half plane $H^2$. Then the product (1) defines a holomorphic function in the half-plane $\mathrm{Re}(s) > 1$. Moreover, $Z_R(s)$ satisfies the functional equation

$$Z_R(s)Z_R(-s) = ((1 - \exp(2\pi is))(1 - \exp(-2\pi is)))^{2-2g}. \tag{3}$$

In particular, for the multiplicity $m_0$ of its singularity in $s = 0$ we have

(4) $$m_0 = 2 - 2g = \chi(X).$$

The functional equation (3) is a consequence of the functional equation

$$Z_S(1-s) = Z_S(s) \exp\left(2(2-2g)\int_0^{s-1/2} (\pi t)\tan(\pi t) dt\right),$$

for the classical Selberg zeta-function

(5) $$Z_S(s) = \prod_c \prod_{N \geq 0} (1 - \exp(-(s+N)l_c)), \qquad \mathrm{Re}(s) > 1, \ s \in \mathbb{C},$$

(see [16]) and the obvious relation $Z_R(s) = Z_S(s+1)/Z_S(s)$.

The Ruelle zeta-function $Z_R$ always satisfies a functional equation similar to (3) relating $Z_R(s)$ to $Z_R(-s)$. However, the general theory of the functional equation for $Z_R$ will not be discussed here. Instead we shall only discuss certain formulas for the multiplicity $m_0$ of the singularity of $Z_R$ in $s = 0$ generalizing formula (4).

To formulate the main result, we shall use group-theoretical descriptions of the geodesic flow and its hyperbolic structure (Anosov property).

Let $g_0$ and $k_0$ be the respective Lie algebras of $G$ and $K$, and let $g_0 = k_0 \oplus p_0$ be the Cartan decomposition of $g_0$ being orthogonal with respect to the Killing form. Identify $p_0$ with the tangent space $T_{eK}(Y)$. Let $a_0 \subset p_0$ be a one-dimensional abelian subspace, and let $M \subset K$ be the centralizer of $a_0$ in $K$.

Now consider the action

$$\Phi_A : A \times \Gamma\backslash G/M \longrightarrow \Gamma\backslash G/M, \quad A = \exp(a_0),$$

defined by

(6) $$(a, \Gamma gM) \mapsto \Phi_a(\Gamma gM) = \Gamma g a^{-1} M, \quad a \in A.$$

The definition of the action $\Phi_A$ is independent of the choice of a metric on $X$, and we shall denote this action in the following as the *abstract* geodesic flow.

Once and for all we fix an orientation of the flow $\Phi_A$ which amounts to fixing an (open) Weyl-chamber $a_0^+$ in $a_0$, called the positive chamber, and we shall restrict attention only to the action of $A^+ = \exp(a_0^+) \subset A$.

Now in terms of the abstract geodesic flow $\Phi_A$ on $\Gamma\backslash G/M$ the hyperbolic structure can be described as follows.

Let $n_0^+$ and $n_0^-$ be the subspaces of $g_0$ on which $ad(X)$ for $X \in a_0^+$ acts by positive and negative eigenvalues $\alpha(X)$, respectively. The nilpotent Lie algebras $n_0^\pm$ are real $MA$-modules with respect to the adjoint action.

Let $\mathcal{P}^\pm$ be the locally homogeneous vector bundles

(7) $$\Gamma\backslash G \times_M (n_0^\pm) \to \Gamma\backslash G/M,$$

and regard these as subbundles of the tangent bundle $T(\Gamma\backslash G/M)$. The real vector bundles $\mathcal{P}^\pm$ are integrable, and the foliation of $\Gamma\backslash G/M$ obtained by integrating $\mathcal{P}^+(\mathcal{P}^-)$ is the $\Phi_A$-invariant unstable (stable) foliation of $\Gamma\backslash G/M$. More precisely, the tangent bundle $T(\Gamma\backslash G/M)$ admits a $d(\Phi_A)$-equivariant decomposition

(8) $$T(\Gamma\backslash G/M) = \mathcal{P}^+ \oplus T^0(\Gamma\backslash G/M) \oplus \mathcal{P}^-$$

into the direct sum of the stable subbundle $\mathcal{P}^-$, the central subbundle $T^0(\Gamma\backslash G/M)$, and the unstable subbundle $\mathcal{P}^+$. Tangent vectors in $\mathcal{P}^-(\mathcal{P}^+)$ are contracted (expanded) exponentially by the differential $d(\Phi_a)$, $a \in A^+$.

Note that the stable and the unstable leaves of $\Gamma\backslash G/M$ are smoothly embedded smooth submanifolds.

Now we shall associate to both foliations of $\Gamma\backslash G/M$ canonical differential forms

$$\Omega_R(\mathcal{P}^\pm) \in C^\infty(\Lambda^{2d-2}T^*(\Gamma\backslash G/M)), \qquad d = \dim(X).$$

We begin with the construction of left $G$-invariant and $\mathrm{End}(\hat{\mathcal{P}}^\pm)$-valued 2-forms

$$\omega_R^\pm \in C^\infty(\Lambda^2 T^*(G/MA) \otimes \mathrm{End}(\hat{\mathcal{P}}^\pm)),$$

where $\hat{\mathcal{P}}^\pm$ denotes the $G$-homogeneous vector bundle

$$G \times_{MA} (n_0^\pm) \to G/MA.$$

By $G$-invariance it suffices to define $\omega_R^\pm$ in $\underline{e} = eMA$.

We choose a (real) basis $\{Z_j\}$ of the space $n_0^\pm$. Let $\{Z^j\}$ be the dual basis of $(n_0^\pm)^*$. Set

(9) $$(\omega_R^\pm)_{\underline{e}} = \left((\omega_R^\pm)_{\underline{e}}(\circ,\circ)_j^k\right)_{j,k}, \quad j, k = 1, \ldots, \dim(n_0^\pm) = d - 1,$$

where

(10) $$(\omega_R^\pm)_{\underline{e}}(\hat{X}_{\underline{e}}, \hat{Y}_{\underline{e}})_j^k = \langle Z^k, -[[X, Y]_0, Z_j]\rangle$$

for $X, Y \in n_0^+ \oplus n_0^-$ and $\hat{X}_{\underline{e}} = d_e(\pi)(X)$ ($\pi: G \to G/MA$ being the canonical projection). Here $[X, Y]_0$ denotes the $m_0 \oplus a_0$-component of $[X, Y]$ with respect to the decomposition $g_0 = (m_0 \oplus a_0) \oplus n_0^+ \oplus n_0^-$, where $m_0$ denotes the Lie algebra of $M$.

Regard the matrix $(\omega_R^\pm)_{\underline{e}}$ as an $\mathrm{End}(n_0^\pm)$-valued alternating 2-form on $T_{\underline{e}}(G/MA)$.

Now extend $(\omega_R^\pm)_{\underline{e}}$ to a $G$-invariant and $\mathrm{End}(\hat{\mathcal{P}}^\pm)$-valued 2-form $\omega_R^\pm$ on $G/MA$.

Next lift $\omega_R^\pm$ via $G/M \to G/MA$ to an $\mathrm{End}(\mathcal{P}^\pm)$-valued $G$-invariant differential 2-form on $G/M$. The latter 2-form drops down to an $\mathrm{End}(\mathcal{P}^\pm)$-valued 2-form on $\Gamma\backslash G/M$ also denoted by $\omega_R^\pm$.

Now define

(11) $$\Omega_R(\mathcal{P}^\pm) = \det\left((i/2\pi)\omega_R^\pm\right).$$

Then the forms $\Omega_R(\mathcal{P}^\pm) \in C^\infty(\Lambda^{2d-2}T^*(\Gamma\backslash G/M))$ are *closed basic* forms with respect to the foliation of $\Gamma\backslash G/M$ by the orbits of the abstract geodesic flow $\Phi_A$, i.e., $\Omega_R(\mathcal{P}^\pm)$ is $\Phi_A$-invariant and $i_{\hat{X}}(\Omega_R(\mathcal{P}^\pm)) = 0$ for all sections $\hat{X}$ of $T^0(\Gamma\backslash G/M)$.

If the dimension of $X$ is even, then it follows that

$$\Omega_R(\mathcal{P}^+) = -\Omega_R(\mathcal{P}^-),$$

and the uniquely determined *real eigenvalue* $\neq 0$ of the $\mathrm{End}(\mathcal{P}^\pm)$-valued 2-form $\omega_R^\pm$ on $\Gamma\backslash G/M$ is a well-defined $\Phi_A$-invariant 2-form

$$\mu_R^\pm \in C^\infty(\Lambda^2 T^*(\Gamma\backslash G/M))$$

on $\Gamma \backslash G/M$. Since $\mu_R^\pm$ is an *exact* 2-form, it follows that there exists a (uniquely determined) left $G$-invariant and right $\Phi_A$-invariant 1-form $\alpha_R^\pm$ on $G/M$ such that $\alpha_R^\pm$ drops down to a $\Phi_A$-invariant 1-form

$$\alpha_R^\pm \in C^\infty(T^*(\Gamma \backslash G/M))$$

which satisfies

(12) $$d\alpha_R^\pm = (i/2\pi)\mu_R^\pm.$$

The main reason to formulate all constructions by using the abstract geodesic flow $\Phi_A$ instead of the geodesic flow $\Phi_t$ is that since the multiplicity $m_0$ is independent of the choice of any scaling of the negative curvature metric on $X$, it is natural to look for a formula for $m_0$ which does not refer to the metric.

Now let $\phi : S(X) \to \Gamma \backslash G/M$ be the diffeomorphism obtained by composing the canonical isomorphism of $S(X)$ and $\Gamma \backslash S(Y)$ with the $G$-equivariant map $S(Y) \ni g(eK, X) \mapsto gM \in G/M$, $X \in a_0^+$ of $S(Y)$ onto $G/M$. The stable and unstable foliations of $\Gamma \backslash G/M$ then obviously correspond to the stable and unstable foliations of $S(X)$ associated to the geodesic flow on $S(X)$.

**Theorem 1.** *Let the dimension of $X$ be even. Then the multiplicity $m_0$ of the singularity of $Z_R$ in $s = 0$ is given by the formula*

(13) $$m_0 = \int_{S(X)} \phi^*(\Omega_R(\mathcal{P}^+) \wedge \alpha_R^-) = \int_{S(X)} \phi^*(\Omega_R(\mathcal{P}^-) \wedge \alpha_R^+).$$

Moreover, the functional equation of $Z_R$ (not given here) implies that *all* singularities of $Z_R$ outside the critical strip $\operatorname{Re}(s) \in [-h, h]$ have multiplicity $2m_0$.

The differential forms $\phi^*(\Omega_R(\mathcal{P}^+) \wedge \alpha_R^-)$ and $\phi^*(\Omega_R(\mathcal{P}^+-) \wedge \alpha_R^+)$ should be regarded as representing a (top-degree) *secondary characteristic class* of the normal bundle of the weak-stable and weak-unstable foliation of $S(X)$, respectively.

There is an equivalent description of $\Omega_R(\mathcal{P}^\pm)$ which emphasizes the analogy of the forms $\Omega_R(\mathcal{P}^\pm)$ with the Pfaffian of the curvature of the Levi-Civita connection of a Riemannian manifold.

In fact, consider the involution $J$ on the tangent bundle $T(G/MA)$ defined by $J|\hat{\mathcal{P}}^\pm = \pm id|\hat{\mathcal{P}}^\pm$ and set $B(X, Y) = \Omega(X, JY)$ for the $G$-invariant symplectic form $\Omega$ on $G/MA$ obtained by reduction of the canonical symplectic form on $T(G/K)\backslash 0$. The leaves of the stable and unstable foliations (of $G/MA$) are Lagrangian submanifolds with respect to $\Omega$. Then $B$ is an invariant pseudo-Riemannian metric of signature $(d-1, d-1)$. The curvature 2-form $\omega_D$ of the corresponding torsion-free pseudo-Riemannian connection $D$ (also considered in [13]) splits as

(14) $$\omega_D = \begin{pmatrix} \omega_D^+ & 0 \\ 0 & \omega_D^- \end{pmatrix}, \qquad \omega_D^- = -\omega_D^+$$

according to the $G$-invariant decomposition $T(G/MA) = \hat{\mathcal{P}}^+ \oplus \hat{\mathcal{P}}^-$. Then the lift (via $G/M \to G/MA$) of the determinant of $(i/2\pi)\omega_D^\pm$ coincides with the form $\Omega_R(\mathcal{P}^\pm)$ up to an exact basic form.

The differential forms $\Omega_R(\mathcal{P}^\pm) \wedge \alpha_R^\pm$ also can be regarded as the top-degree components of $A$-equivariantly closed forms (of mixed degree) on $\Gamma \backslash G/M$. Therefore Theorem 1 can be regarded as a *regularized* analog of a localization formula in equivariant cohomology (see [1]). In particular, it is natural to regard the multiplicity formula (13) as an analog of the Poincaré-Hopf formula for the sum of the indices of a (non-degenerate) vector field.

In the case of a compact Riemannian surface $X = \Gamma \backslash H^2$ equation (4) follows from Theorem 1 by an elementary calculation. In fact, in this case the differential forms

$$(4\pi^2)\phi^*(\Omega_R(\mathcal{P}^+) \wedge \alpha_R^-)$$

represent the Godbillon-Vey class of the weak-stable foliation of the geodesic flow on the sphere bundle $S(X)$ of the surface $X$. Thus Theorem 1 can be interpreted as the assertion that the Godbillon-Vey invariant of the weak-stable foliation can be calculated from the closed orbits of the geodesic flow. The relation between the Godbillon-Vey invariant of the weak-stable foliation of $S(X)$ and the Euler characteristic of the surface $X$ is a well-known famous result due to Roussarie (see [9]).

By the method of symbolic dynamics it follows that $m_0$ coincides with the alternating sum of the (finite) dimensions of the generalized eigenspaces of the transfer operators $L_*(0)$ for the eigenvalue 1. Recall that the operators $L_*(0)$ do *not* depend on return times (see [5])! Therefore Theorem 1 also can be regarded as a formula for the integer associated to the flow by forming this *analytical index*. Although the definition of the latter index strongly depends on the choice of a Markov-family of local sections, it is, in fact, independent of the ambiguities involved in the construction of the local sections. Moreover, Theorem 1 shows that the analytical index coincides with the integer defined by integrating a secondary characteristic cohomology class that depends only on the hyperbolic structure of the flow.

Theorem 1 is but a special case of more general formulas relating the multiplicities of the singularities of generalized Selberg zeta-functions at special points to integrals of canonically associated secondary characteristic classes.

Next we combine Theorem 1 with proportionality theory. Let $Y^d$ be the compact dual symmetric space of $Y$. $Y^d$ is a rank one space, and all geodesics are closed and have the same length. Let $Y^d_{\text{geo}}$ be the space of all (oriented) geodesics in $Y^d$ (see [2]).

**Theorem 2.** *Let $X$ be as in Theorem 1. Then*

$$(15) \qquad m_0 = \left(\chi(X)/\chi(Y)^d\right) \chi(Y^d_{\text{geo}}),$$

*where $\chi$ always denotes Euler characteristic.*

Moreover, by a calculation of the Euler characteristics of the compact homogeneous spaces occurring in (15) one obtains

**Corollary 1.** *Let $X$ be as in Theorem 2. Then*

$$(16) \qquad m_0 = (\dim(X)/2)\chi(X).$$

Note that, in contrast to the even-dimensional case, for an odd-dimensional *real* hyperbolic space we have the following formula.

**Theorem 3.** *Let $X = \Gamma \backslash H^{2n+1}$ be a compact real hyperbolic space of dimension $2n+1$. Then the multiplicity of the singularity of $Z_R$ in $s = 0$ is given by*

$$(17) \qquad 2\left((-1)^{n+1}b_{n+1}(X) + \cdots + (-1)^{2n+1}(n+1)b_{2n+1}(X)\right)$$

where $b_p(X)$ is the $p$th Betti-number of $X$.

*On the proofs.* Our proof of Theorem 1 is, unfortunately, much more complicated than the result itself suggests. Thus we can give here only some hints on how the assertion can be proved. More details can be found in [12].

The proof rests on a cohomological trace formula which can be regarded as a common (non-commutative) generalization of the Poisson summation formula and the Lefschetz fixed-point formula. It implies that the zeta-function $Z_R$ is closely related to the alternating product of infinite-dimensional (regularized) characteristic determinants of global Frobenius-operators (canonically determined by the action of the geodesic flow) on the cohomology groups of some differential complexes associated to the invariant foliations of $S(X)$. More precisely, it yields a cohomological formula for the multiplicity $m_0$ in terms of the Lie-algebra-cohomology of $n^\pm$ with values in the Harish-Chandra modules of the irreducible representations of $G$ in $L^2(\Gamma\backslash G)$.

Now if $\dim(X)$ is even, then $G$ has a compact Cartan subgroup $H$. The cohomological formula for $m_0$ turns out to be connected with an analogous cohomological formula obtained by replacing $n^\pm$-cohomology by Lie-algebra-cohomology with respect to the nilradicals of the (complex) Borel algebras containing the complexified Lie algebra of $H$. This can be proved by using Osborne's character-formula (see [10]) and suitable patching conditions for characters on neighbouring Cartan subgroups. But the latter number coincides with the analytical index of the (elliptic) deRham complex on the space $\Gamma\backslash G/H$. By working backwards with the corresponding index-form (given by Gauss-Bonnet), one finally ends up with the formula (13).

If the dimension of $X$ is odd, then one can explicate the cohomological formula for $m_0$ directly by using results from [4].

It is also possible to give more traditional proofs by applying suitable explicit Selberg trace formulas. A proof of Theorem 3 along these lines can be found in [6]. In the even-dimensional case a proof of Corollary 1 resting on Selberg trace formulas is given in [3]. Finally, for quaternionic hyperbolic spaces a representation theoretical proof of Corollary 1 is given in [15].

## References

[1] N. Berline, E. Getzler, and M. Vergne, *Heat kernel and Dirac operators*, Springer, Berlin, 1992.

[2] A. Besse, *Manifolds all of whose geodesics are closed*, Ergeb. Math. Grenzgeb. (3), vol. 93, Springer, Berlin, 1978.

[3] U. Bunke and M. Olbrich, *Theta and zeta functions for locally symmetric spaces of rank one*, preprint SFB 288, Berlin, 1994.

[4] D. Collingwood, *Representations of rank one Lie groups*, Pitman Res. Notes Math. Ser., vol. 137, Pitman, Boston, 1985.

[5] D. Fried, *The zeta functions of Ruelle and Selberg. I*, Ann. Sci. École Norm. Sup. (4) **19** (1986), 491–517.

[6] \_\_\_\_\_\_, *Analytic torsion and closed geodesics on hyperbolic manifolds*, Invent. Math. **84** (1986), 523–540.

[7] R. Gangolli, *Zeta functions of Selberg's type for compact space forms of symmetric spaces of rank one*, Illinois J. Math. **21** (1977), 1–41.

[8] C. Godbillon and J. Vey, *Un invariant de feuilletages de codimension* 1, C. R. Acad. Sci. Paris **273** (1971), 92–95.

[9] S. Hurder and A. Katok, *Differentiability, rigidity and Godbillon-Vey classes for Anosov flows*, Inst. Hautes Études Sci. Publ. Math. (1990), 5–61.

[10] H. Hecht and W. Schmid, *Characters, asymptotics and n-homology of Harish-Chandra modules*, Acta Math. **151** (1983), 49–151.

[11] A. Juhl, *On the functional equation of dynamical theta functions*, preprint IAAS, Berlin, 1993.

[12] ______, *Zeta-Funktionen, Indextheorie und hyperbolische Dynamik*, Thesis, Humboldt-Universität, Berlin, 1993.

[13] M. Kanai, *Geodesic flows of negatively curved manifolds with smooth stable and unstable foliations*, Ergodic Theory Dynamical Systems **8** (1988), 215–239.

[14] D. Ruelle, *Zeta functions for expanding maps and Anosov flows*, Invent. Math. **34** (1976), 231–242.

[15] S. Seifarth, *Kohomologische Untersuchungen dynamischer Zeta-Funktionen*, Thesis, Humboldt-Universität, Berlin, 1994.

[16] A. Selberg, *Harmonic analysis and discontinuous groups in weakly symmetric Riemannian spaces with applications to Dirichlet series*, J. Indian Math. Soc. **20** (1956), 47–87.

[17] M. Wakayama, *Zeta functions of Selberg's type associated with homogeneous vector bundles*, Hiroshima Math. J. **15** (1985), 235–295.

WIAS, Mohrenstrasse 39, 10117 Berlin, Germany
*E-mail address*: juhl@iaas-berlin.d400.de